\documentclass[12pt]{amsart}
\usepackage{amssymb,amsthm}
\usepackage{enumitem}

\theoremstyle{plain}

\newtheorem*{lem*}{Lemma}
\newtheorem*{prop*}{Proposition}
\newtheorem*{thm*}{Theorem}
\newtheorem*{cor*}{Corollary}
\newtheorem*{conj*}{Conjecture}
\theoremstyle{remark}

\newtheorem*{rmks}{Remarks}

\newcommand{\Z}{{\mathbb{Z}}}

\renewcommand{\le}{\leqslant}
\renewcommand{\ge}{\geqslant}
\newcommand{\GL}{\mathrm{GL}}
\newcommand{\SL}{\mathrm{SL}}

\newcommand{\F}{\mathcal{F}}
\newcommand{\St}{\mathrm{St}}

\begin{document}
\title[Factoring tilting modules]{Factoring tilting modules for
  algebraic groups} \author{S.R.~Doty} \address{Mathematics and
  Statistics, Loyola University Chicago, Chicago, IL 60626, USA}
\subjclass[2000]{} \date{29 June 2009}
\thanks{Research supported by the Mercator Programme, DFG}

\begin{abstract}
Let $G$ be a semisimple, simply-connected algebraic group over an
algebraically closed field of characteristic $p>0$.  We observe that
the tensor product of the Steinberg module with a minuscule module is
always indecomposable tilting. Although quite easy to prove, this fact
does not seem to have been observed before. It has the following
consequence: If $p \ge 2h-2$ and a given tilting module has highest
weight $p$-adically close to the $r$th Steinberg weight, then the
tilting module is isomorphic to a tensor product of two simple
modules, usually in many ways.
\end{abstract}
\maketitle

\noindent
Let $G$ be a semisimple, simply-connected algebraic group over an
algebraically closed field $k$ of characteristic $p>0$. For
convenience we assume the underlying root system is indecomposable.
Tensor products are over $k$ unless otherwise specified. Fix a maximal
torus $T$ in $G$ and write $X(T)$ for the character group of $T$. Note
that $X(T) \simeq \Z^n$ for some $n$.  By ``$G$-module'' we mean
``rational $G$-module''. Fix a Borel subgroup $B$ containing $T$ and
let the negative roots be determined by $B$. Let $$X(T)^+ = \{ \lambda
\in X(T): (\alpha^\vee, \lambda) \ge 0, \text{ all simple roots }
\alpha \}$$ be the set of dominant weights and $$X_r(T) = \{\lambda
\in X(T)^+: (\alpha^\vee, \lambda) < p^r, \text{ all simple roots }
\alpha \}.$$ The set $X_1(T)$ is known as the restricted region and
its elements are often called restricted weights. For any $\lambda \in
X(T)^+$ let

$\Delta(\lambda) = $ the Weyl module of highest weight $\lambda$;

$\nabla(\lambda) = $ the dual Weyl module of highest weight $\lambda$;

$L(\lambda) = $ the simple $G$-module of highest weight $\lambda$.

\noindent
The main properties of these families of modules are summarized in
\cite{Jantzen}, to which the reader should also refer for any
unexplained notation or terminology.  

Let $\F(\Delta)$ be the full subcategory of the category of
$G$-modules whose objects have an ascending filtration with successive
sub-quotients isomorphic to various Weyl modules; $\F(\nabla)$ is
defined similarly with $\nabla$ in place of $\Delta$.  Recall that the
objects of $\F(\Delta) \cap \F(\nabla)$ are called tilting modules and
the category of tilting modules is closed under tensor products,
direct sums, and direct summands. For each $\lambda \in X^+$, there is
a unique (up to isomorphism) indecomposable tilting module of highest
weight $\lambda$, denoted by $T(\lambda)$. Every tilting module is
isomorphic to a direct sum of various $T(\lambda)$. Since
$\Delta(\lambda)$ is isomorphic to the contravariant dual of
$\nabla(\lambda)$ it follows immediately that whenever
$\Delta(\lambda)$ is simple as a $G$-module, then
\begin{equation}\label{eq:1}
L(\lambda) \simeq \Delta(\lambda) \simeq \nabla(\lambda) \simeq
T(\lambda).
\end{equation}
Conversely, any simple tilting module must be a simple Weyl module.

A dominant weight is called \emph{minuscule} if the weights of
$\Delta(\lambda)$ form a single orbit under the action of the Weyl
group $W$. This forces $\Delta(\lambda)$ to be simple, so \eqref{eq:1}
holds for any minuscule weight $\lambda$. When $\lambda$ is minuscule
we shall refer to any of the isomorphic modules in \eqref{eq:1} as a
minuscule module. Note that the zero weight is minuscule and the
trivial module is a minuscule module by our definition.  Minuscule
weights are classified in \cite[ch.~VIII, prop.~7]{Bourbaki}.  For the
reader's convenience we list them in Table \ref{tab:1}.
\begin{table}[ht]
\begin{tabular}{|llll|} \hline
  Type & Highest Weight & Dimension & Name \\ \hline
  $A_n$ & $\varepsilon_j$ $(1 \le j \le n)$ & $\binom{n+1}{j}$ & 
        exterior powers of natural \\
  $B_n$ & $\varepsilon_n$ & $2^n$ & spin \\
  $C_n$ & $\varepsilon_1$ & $2n$ & natural \\
  $D_n$ & $\varepsilon_1, \varepsilon_{n-1}, \varepsilon_n$  
        & $2n$, $2^{n-1}$, $2^{n-1}$ 
        & natural, $\frac{1}{2}$-spin, $\frac{1}{2}$-spin \\
  $E_6$ & $\varepsilon_1, \varepsilon_6$ & 27, 27 & minimal \\
  $E_7$ & $\varepsilon_7$ & 56 & minimal \\
  $E_8$ & none &&\\
  $F_4$ & none &&\\
  $G_2$ & none &&\\ \hline
\end{tabular}
\caption{Minuscule modules} \label{tab:1}
\end{table}
In the table, $\varepsilon_1, \dots, \varepsilon_n$ are the
fundamental weights, defined by the requirement $(\alpha_i^\vee,
\varepsilon_j) = \delta_{i,j}$ for all $i,j$ (with respect to the
usual ordering of the simple roots). Note that all minuscule weights
belong to the restricted region $X_1(T)$ for any $p$.

Let $\rho \in X(T)$ be half the sum of the positive roots. Write
$\St_r := \Delta((p^r-1)\rho)$ for the $r$th Steinberg module; this is
a simple tilting module for every $r>0$. We write $\St$ for $\St_1$.

\begin{lem*}
  If $\lambda$ is minuscule then
  $\St \otimes L(\lambda) \simeq T((p-1)\rho + \lambda).$
\end{lem*}

\begin{proof}
  In \cite[Proposition 5.5]{Donkin:TiltHandbook} it is proved (by an
  application of Brauer's formula) that if $(\alpha_0^\vee, \lambda)
  \le p$, where $\alpha_0$ is the highest short root, then the
  character of $T((p-1)\rho+\lambda)$ is equal to the character of
  $\St$ multiplied by the character of the orbit of $\lambda$ under
  the action of $W$. Now the tensor product $\St \otimes L(\lambda)$
  in question is the tensor product of two tilting modules, hence is
  itself tilting. By highest weight considerations a copy of
  $T((p-1)\rho + \lambda)$ must occur as a direct summand.  Thus we
  are done once we have verified that $(\alpha_0^\vee, \lambda) \le
  p$. But this is easy to check, by comparing the classification of
  minuscule weights in Table \ref{tab:1} with a list of highest short
  roots (see \cite[\S12, Table 2]{Humphreys}).
\end{proof}

We now want to generalize the above result. Say that a weight
$\lambda$ is $r$-minuscule if $\lambda$ can be written in the form
$\lambda = \sum_{j=0}^{r-1} p^j \lambda^j$, where each $\lambda^j$ is
minuscule. For such $\lambda$ we obviously have 
\[
  L(\lambda) \simeq L(\lambda^0) \otimes L(\lambda^1)^{[1]} \otimes
  \cdots \otimes L(\lambda^{r-1})^{[r-1]}
\]
by Steinberg's tensor product theorem. 

Let $h$ be the Coxeter number of the underlying root system.  Recall
(Donkin \cite[p.~47, Example 1]{Donkin:Zeit}) that if $p \ge 2h-2$ and
$\lambda \in X_r(T)$ then $T((p^r-1)\rho + \lambda)$ is isomorphic to
the projective cover of $L((p^r-1)\rho + w_0\lambda)$ in the category
of $G_r$-modules. Here $G_r$ is the $r$th Frobenius kernel of $G$ and
$w_0$ is the longest element of the Weyl group.  Donkin has
conjectured that this holds for any $p$; see
\cite[(2.2)]{Donkin:Zeit}. He proved in \cite[(2.1)]{Donkin:Zeit} (see
also \cite[II.E.9]{Jantzen}) that $T(\tau) \otimes T(\mu)^{[r]}$ is
tilting, for any $\tau \in (p^r-1)\rho + X_r(T)$, $\mu \in X(T)^+$,
and morover if $p \ge 2h-2$ (or if the conjecture holds for $p <
2h-2$) then
\begin{equation}\label{eq:2}
  T(\tau) \otimes T(\mu)^{[r]} \simeq T(\tau + p^r \mu).
\end{equation}
This statement is known as the tensor product theorem for tilting
modules.

\begin{prop*} 
  Assume that Donkin's conjecture holds for $G$ if $p < 2h-2$. If
  $\lambda$ is $r$-minuscule and $\mu\in X(T)^+$
  then $$T(\mu)^{[r]}\otimes \St_r \otimes L(\lambda) \simeq T(p^r \mu
  + (p^r-1)\rho + \lambda).$$
\end{prop*}

\begin{proof}
  By Steinberg's tensor product theorem it follows that 
  \[
   \St_r \otimes L(\lambda) \simeq \textstyle \bigotimes_{j=1}^{r-1}
   \big(\St \otimes L(\lambda^j)\big)^{[j]} 
  \]
  where $\lambda = \sum_j \lambda^j p^j$ (with $\lambda^j \in X_1(T)$
  for all $j$) is the $p$-adic expansion of $\lambda$. By the lemma we
  get
  \[
   \St_r \otimes L(\lambda) \simeq \textstyle \bigotimes_{j=1}^{r-1}
   \big(T((p-1)\rho+\lambda^j)\big)^{[j]} 
  \]
  and by the tensor product theorem for tilting modules (see
  \eqref{eq:2}) applied inductively it follows that
  \[
   \St_r \otimes L(\lambda) \simeq T((p^r-1)\rho+\lambda).
  \]
  Now tensor both sides by $T(\mu)^{[r]}$ and apply the tensor product
  theorem for tilting modules again to obtain the result.
\end{proof}

In general one would like to understand the indecomposable direct
summands of modules of the form $L \otimes M$ where $L$ is simple and
$M$ is either simple or tilting. The proposition provides many
examples where such tensor products are in fact indecomposable tilting
modules.

\begin{cor*}
  Assume that Donkin's conjecture holds for $G$ if $p < 2h-2$. If
  $\lambda$ is $r$-minuscule and $\mu\in X(T)^+$ then:
  \begin{enumerate}[label={\rm(\alph*)},leftmargin=*,itemsep=0.5em]
  \item $T(p^r\mu + (p^r-1)\rho) \otimes L(\lambda) \simeq T(p^r \mu +
    (p^r-1)\rho + \lambda)$.

  \item If $T(\mu)$ is simple then $\St_r \otimes L(p^r\mu+\lambda)
    \simeq T(p^r \mu + (p^r-1)\rho + \lambda)$.
  \end{enumerate}
\end{cor*}

\begin{proof}
  By the tensor product theorem for tilting modules we have
  $T(\mu)^{[r]} \otimes \St_r \simeq T(p^r\mu + (p^r-1)\rho)$. This
  proves (a).

  If $T(\mu) \simeq L(\mu)$ then $T(\mu)^{r]} \otimes L(\lambda)
    \simeq L(\mu)^{r]} \otimes L(\lambda) \simeq L(p^r\mu + \lambda)$,
  by Steinberg's tensor product theorem. This proves (b).
\end{proof}

\begin{rmks}
1. In case $G$ is of Type $A_1$ or $A_2$ it is known that Donkin's
conjecture holds for all $p$.

2. Given two simple modules $L,M$ one may express each one as a twisted
tensor product of restricted simple modules
\begin{align*}
L &\simeq L_0\otimes L_1^{[1]} \otimes L_2^{[2]} \otimes \cdots \\
M &\simeq M_0\otimes M_1^{[1]} \otimes M_2^{[2]} \otimes \cdots 
\end{align*}
by Steinberg's tensor product theorem. Interchanging $L_j$ and $M_j$
in arbitrary selected positions $j$ results in two new simple modules
$L', M'$ such that $L \otimes M \simeq L' \otimes M'$. This is
immediate by commutativity of tensor product.  Applying this
observation to the pair $\St_r$, $L(p^r\mu+\lambda)$ in part (b) of
the corollary produces many factorizations $$T(p^r\mu +
(p^r-1)\rho+\lambda) \simeq L(\lambda') \otimes L(\mu')$$ where
$\lambda', \mu'$ are the highest weights of the rearranged tensor
products.

3. There exist factorizations of indecomposable tilting modules not of
the form in the proposition or corollary. For example, for $G = \SL_3$
in characteristic 3 one has from \cite[5.2]{DM} the factorization
$T(3,0) \simeq L(2,0) \otimes L(1,0)$.

4. The proposition and corollary are most effective when $p$ is small,
in which case more weights are close to a Steinberg weight. For
instance, for $G = \SL_2$ in characteristic $2$ it follows from the
corollary that every indecomposable tilting module is isomorphic to a
tensor product of two simple modules. This was known previously; see
\cite[2.5]{DH}.

5. The above results can be extended to the reductive case by the
usual arguments, although if $p=2$ one has to deal with the
possibility that the Steinberg module (as defined above) may fail to
exist; e.g. consider $G=\GL_2$. Thus one may need to pass to a
covering. We leave the details to the reader. 

Another possibility in the reductive case is to replace $\rho$ by any
weight $\rho'$ satisfying the condition $(\alpha^\vee, \rho') = 1$ for
all simple roots $\alpha$. The resulting modules
$\Delta((p^r-1)\rho')$ satisfy the desired properties of Steinberg
modules and may be used in place of the $\St_r$ in the above
arguments; see the remark in \cite[II.3.18]{Jantzen}. In case
$G=\GL_n$ one may wish to apply this remark with the weight $\rho' =
\sum_{j=1}^n (n-j)\varepsilon_j$ replacing $\rho$.

6. The tensor product $\St \otimes L(\lambda)$ considered in the lemma
(for minuscule $\lambda$) is projective as a $kG^F$-module, where $F$
is the $p$-Frobenius endomorphism of $G$ and $G^F$ is the finite
Chevalley group of $F$-fixed points in $G$.  Moreover, if $U(\mu)$
denotes the projective cover of $L(\mu)$ in the category of
$kG^F$-modules (for $\mu \in X_1(T)$) then $U((p-1)\rho + w_0\lambda)$
occurs once (see \cite{Jantzen2}, \cite{Chast}) as a direct summand of
$\St \otimes L(\lambda)$, viewed as a $kG^F$-module. It is natural to
ask if $U((p-1)\rho + w_0\lambda) \simeq \St \otimes L(\lambda)$. This
is not always the case; see \cite{Tsushima} for a study of this
question.
\end{rmks}

\end{document}